\documentclass[12pt, reqno]{amsart}

\usepackage{amssymb, amsmath, amsthm}
\usepackage[backref]{hyperref}
\usepackage[alphabetic,backrefs,lite]{amsrefs}
\usepackage{amscd}   
\usepackage[all]{xy} 

\DeclareFontEncoding{OT2}{}{} 

\usepackage[usenames,dvipsnames]{color}

\theoremstyle{plain}
\newtheorem{prop}{Proposition}[section]

\newtheorem{thm}[prop]{Theorem}
\newtheorem{cor}[prop]{Corollary}

\newtheorem{lemma}[prop]{Lemma}

\theoremstyle{definition}

\newtheorem{rmk}[prop]{Remark}

\newtheorem{example}[prop]{Example}

\newcommand{\PP}{{\mathbb P}}
\newcommand{\C}{{\mathbb C}}

\newcommand{\R}{{\mathbb R}}
\newcommand{\Z}{{\mathbb Z}}

\newcommand{\frakd}{{\mathfrak d}}

\newcommand{\frakB}{{\mathfrak B}}

\newcommand{\frakG}{{\mathfrak G}}

\newcommand{\frakI}{{\mathfrak I}}
\newcommand{\frakJ}{{\mathfrak J}}
\newcommand{\frakK}{{\mathfrak K}}
\newcommand{\frakL}{{\mathfrak L}}

\newcommand{\frakQ}{{\mathfrak Q}}

\newcommand{\frakT}{{\mathfrak T}}

\DeclareMathOperator{\rk}{rk}

\DeclareMathOperator{\Aut}{Aut}

\DeclareMathOperator{\Pic}{Pic}

\DeclareMathOperator{\Bir}{Bir}
\DeclareMathOperator{\NS}{NS}

\numberwithin{equation}{section}
\numberwithin{table}{section}

\title[Degree two unirational parametrizations]{Degree two unirational parametrizations over the real field}
\author{Brendan Hassett}
\author{Hayato Takagi}
\author{Sho Tanimoto}

\address{Department of Mathematics, Brown University Box 1917, 151 Thayer Street, Providence, RI 02912, USA}
\email{brendan\_hassett@brown.edu}

\address{Graduate School of Mathematics, Nagoya University, Furocho Chikusa-ku, Nagoya, 464-8602, Japan}
\email{hayato.takagi.c6@math.nagoya-u.ac.jp}

\address{Graduate School of Mathematics, Nagoya University, Furocho Chikusa-ku, Nagoya, 464-8602, Japan}
\email{sho.tanimoto@math.nagoya-u.ac.jp}

\begin{document}


	\maketitle
	
	
	\section{Introduction}

    An algebraic variety is rational if it is birational to projective space
    and unirational if it is dominated
    by a rational variety. There exist complex threefolds that are unirational but not rational (\cites{IM71, CG72, AM72}) A smooth projective complex variety is rationally connected if any two points may be connected by a rational curve. One outstanding question is whether there exists a rationally connected smooth projective variety which is not unirational. In dimensions $\leq 2$, those three notions are equivalent over $\mathbb C$, but it is expected that unirationality is a stronger notion than rational connectedness in higher dimensions.

On the other hand, a variety defined over $\R$ may be rational over $\mathbb C$ but not over $\R$, e.g., a conic $Q_{3, 0} = \{x^2 + y^2 + z^2 =0\} \subset \PP^2$ is isomorphic to $\PP^1$ over $\mathbb C$ but has no real point, so it cannot be rational over $\R$. More generally, if a smooth projective variety $X$ is unirational over $k$ then $X(k)$ is non-empty. Another obstruction comes from Galois cohomology: 
If we have a unirational parametrization of degree $n$ then $n$ must annihilate 
$H^1(\mathrm{Gal}(\overline{k}/k), \Pic (X_{\overline{k}}))$
as for any degree-$n$ morphism $f : Y \to X$, we have $f_*f^* = [n]$ on $\Pic (X_{\overline{k}})$.
When $k = \R$, the Galois cohomology $H^1(\mathrm{Gal}(\mathbb C/\R), \Pic (X_{\mathbb C}))$ is a $2$-torsion group. Indeed, this follows from a fact that $$H^1(\mathbb Z/2\mathbb Z = \langle \tau\rangle, M) = \{ m \in M \, | \, \tau m = -m\}/\{\tau m - m \, |\, m \in M\}.$$
Thus one may ask if there exists a geometrically unirational variety, defined over $\R$ with real
points, which does not admit a degree-two real unirational parametrization. 

A systematic study of such questions has been conducted by Trepalin in \cite{Tre19} under the name of {\it Galois unirational surfaces}, and he obtained a complete answer for surfaces with rational points over fields of characteristic $0$. In this survey paper, we illustrate the classification of Galois unirational surfaces.
The following theorem is a corollary of \cite{Tre19}*{Theorem 1.15}, and we give a different proof:

\begin{thm}[Trepalin]
\label{thm:intromain}
    Let $S$ be a smooth real projective surface such that $S_\C$ is rational.
    Then there exists a degree-two dominant rational map $\PP^2 \dashrightarrow S$ if and only if $S$ is either rational or birational to a minimal conic bundle with a real point.
\end{thm}

When $S$ is a smooth projective geometrically rational surface over $\R$, its birational model is completely classified by the following theorem:

\begin{thm}[{\cite{Kollar}, Corollary 3.4, and Theorems 6.3 and 6.8}]
\label{thm:kollar}
Let $S$ be a smooth real projective surface such that $S_\C$ is rational.
Then $S$ is birationally equivalent over $\R$ to a surface in exactly one of the following classes:
\begin{enumerate}
\item[(1)] $Q_{3,0} \times \PP^1$.
In this case $S(\R)=\emptyset$.
\item[(2)] $\PP^2$.
In this case $S(\R)$ is connected.
\item[$(3_m)$] Minimal conic bundle with $2m\ (m \geq 2)$ singular fibers.
In this case $S(\R)$ has $m$ connected components.
\item[(4)] Minimal del Pezzo surface of degree $2$ and Picard rank $1$.
In this case $S(\R)$ has $4$ connected components.
\item[(5)] Minimal del Pezzo surface of degree $1$ and Picard rank $1$.
In this case $S(\R)$ has $5$ connected components.
\end{enumerate}
\end{thm}

\begin{rmk}
    Minimal conic bundles only admit an even number of singular fibers, and any minimal conic bundle with $2$ singular fibers is rational. 
\end{rmk}

In particular Theorems~\ref{thm:intromain} and \ref{thm:kollar} show the following corollary 
\begin{cor}
\label{cor:degree12}
    Any minimal del Pezzo surface of degree $1$ or $2$ 
    of Picard rank $1$ over $\R$ does not admit a degree-two real unirational parametrization.
\end{cor}

Another corollary of Theorems~\ref{thm:intromain} is the following:
\begin{cor}
\label{cor:degreegeq3}
    Let $S$ be a real del Pezzo surface of degree $\geq 3$.
    If $S(\R) \neq \emptyset$ then $S$ admits a degree-two real unirational parametrization.
\end{cor}

Our proof of Theorems~\ref{thm:intromain} is based on the following idea: Suppose that we have a degree-two unirational parametrization $\PP^2 \dashrightarrow S$. After resolving the indeterminacy and taking the Stein factorization, we obtain a degree-two finite morphism $X' \to S$ where $X'$ is a normal rational projective surface. This induces an involution $\tau'$ on $X'$. After resolving singularities equivariantly and applying the equivariant minimal model program, we obtain an equivariant minimal model $(X, \tau)$ such that $X$ is a smooth rational projective surface. Such pairs are classified by \cite{CMYZ24} which studied conjugacy classes of involutions in the real plane Cremona group. Thus we go through the classification by \cite{CMYZ24} and confirm that in all cases the quotient $X/\tau$ is either rational or birational to a minimal conic bundle. It is obvious that any rational surface admits a degree-two unirational parametrization. For conic bundles, we use an argument by Koll\'ar and Szab\'o. See Corollary~\ref{cor:conicbundleunirational}.

\subsection*{Related results}
We mention results related to this paper.
\cite{CT19}*{Theorem 1.2} shows that if any real geometrically rational surface admits an odd-degree unirational parametrization, then it must be rational over the ground field. \cite{KM17} proved unirationality for del Pezzo surfaces of degree $1$ with conic bundle structures. Knecht used a similar idea to this paper in \cite{Knecht} and proved that any minimal cubic surface over a field of odd characteristic does not admit a degree-two unirational parametrization. For recent results
on rationality of higher-dimensional varieties over $\R$ and other non-closed fields see \cites{HT21b, HT21, HKT22, BW20, BW23, KP23, FJSVV}, for example.

\subsection*{Structure of the paper}
In Section~\ref{sec:equivariantmmp}, we discuss equivariant minimal model program for involutions over non-closed fields. This is discussed over an algebraically closed field in \cite{BBAsian}, and we give a survey for an analogous theory over non-closed fields and the real field. Strictly speaking this section is not necessary for the rest of the paper, and people who are familiar with the equivariant minimal model program may skip this section. In Section~\ref{sec:Kollar} we give an argument by Koll\'ar and Szab\'o and prove that any real conic bundle with a real point admits a degree-two unirational parametrization.
In Section~\ref{sec:main} we prove Theorem~\ref{thm:intromain}.

    \bigskip

\
	
	\noindent
{\bf Acknowledgments:}
This paper stems from the working seminar on rational surfaces over non-closed fields which was held at Rice University during Spring 2015. The authors would like to thank the participants of that seminar.
We would like to thank Andrey Trepalin for pointing
out relevant results in \cite{Tre19} and \cite{Isk67}.
Finally we would like to thank the referee for a careful reading of the manuscript and suggestions to improve the exposition of the paper.

Brendan Hassett was partially supported by the 
Simons Foundation Award 546235 and 
the US National Science Foundation Grant 1929284.
Hayato Takagi was partially supported by JST SPRING, Grant number JPMJSP2125.
Sho Tanimoto was partially supported by JST FOREST program Grant number JPMJFR212Z, by JSPS KAKENHI Grant-in-Aid (B) 23H01067, by JSPS Bilateral Joint Research Projects Grant number JPJSBP120219935, Lars Hesselholt's Niels Bohr professorship, by MEXT Japan, Leading Initiative for Excellent Young Researchers (LEADER), by Inamori Foundation, and by JSPS KAKENHI Early-Career Scientists Grant numbers 19K14512.

\section{The classification of surfaces with biregular involutions over non-closed fields}
\label{sec:equivariantmmp}

The classification of rational surfaces with biregular involutions has been developed over an algebraically closed field of characteristic $\neq 2$ in \cite{BBAsian}. Here we survey the minimal model program of smooth projective surfaces with biregular involutions over a non-closed field, following \cite{BBAsian}. Suppose that $k$ is a perfect field of characteristic $\neq 2$. We denote its absolute Galois group by $G = \mathrm{Gal} (\overline{k}/k)$.

Let $X$ be a smooth projective surface over $k$. 
Let $\NS(X)$ denote the subgroup of the N\'eron-Severi group
$\NS(\overline{X})$ corresponding to divisors defined over $k$,
a finite-index subgroup of $\NS(\overline{X})^G$.
Let $\overline{\mathrm{NE}}$ denote the closed cone of 
effective curves in $\NS(\overline{X})\otimes {\R}$.

Suppose that we have a non-trivial biregular involution $\tau : X \rightarrow X$. A pair $(X, \tau)$ is minimal if any birational morphism $f: X\rightarrow X'$, with a biregular involution $\tau' : X' \rightarrow X'$ such that $\tau' \circ f = f \circ \tau$,  is an isomorphism.

\begin{lemma}
A pair $(X, \tau)$ is minimal if and only if there are no exceptional curves $E$, e.g., a $(-1)$-curve, whose Galois conjugates are disjoint to each other, and the sum of all Galois conjugates $L = \sum g_i\cdot E$ satisfies either $L = \tau L$ or $\tau L \cap L = \emptyset$.
\end{lemma}

\begin{proof}
Suppose that there is an exceptional curve $E$ whose Galois conjugates are disjoint and we have $\tau L = L$ or $\tau L \cap L = \emptyset$. Then we can contract the Galois orbit of $E \cup \tau E$ to $X'$ and moreover, the involution $\tau$ induces an involution $\tau'$ on $X'$. Hence $X$ is not minimal.

Assume that $X$ is not minimal, i.e., it admits a non-trivial birational morphism $f :(X, \tau) \rightarrow (X', \tau')$. Let $E$ be an exceptional curve contracted by $f$. Then its Galois conjugates are also contracted by $f$, so they must be disjoint. Moreover $\tau L + L$ also has to be contracted, so we have $(\tau L + L)^2 <0$. This implies that $\tau L \cdot L \leq0$. Hence $\tau L$ and $L$ coincide or they do not meet.
\end{proof}

\subsection*{General fields}
Our goal here is to give a rough classification of minimal pairs $(X, \tau)$. First we consider the case of $\mathrm{rk} \, \NS(X)^{\tau} >1$:

\begin{lemma}
\label{lemma:bpfpencil}
Let $(X, \tau)$ be a minimal pair with $\mathrm{rk} \, \NS(X)^{\tau} >1$ and $K_X$ not nef. Then $X$ admits a base point free linear system of Iitaka dimension $1$, stable under $\tau$ and defined over $k$.
\end{lemma}

\begin{proof}
Suppose that $\mathrm{rk} \,\NS(X)^{\tau} >1$. Since the canonical class $K_X$ is not nef, the cone theorem implies that
\[
\overline{\mathrm{NE}}(\overline{X})^G = \overline{\mathrm{NE}}(\overline{X})^G_{K_X \geq0} + \sum_{L \in \mathcal E} \R_+ [L]
\]
where $\mathcal E$ is a non-empty countable set and $L$ is the sum of all Galois conjugates of a curve $E$ where $E$ generates a $K_X$-negative ray for $\overline{\mathrm{NE}}(\overline{X})$. It follows from the proof of \cite{HassettRational}*{Theorem 3.9} that $L$ has the following possibilities: (i) $E$ is a $(-1)$-curve and its Galois conjugates are disjoint; (ii) $E$ is a $(-1)$-curve, the Galois orbit of $E$ decomposes as
\[
\{E_1, E_1'\}, \cdots, \{E_r, E_r'\}
\]
where $E_i \cdot E_i' = 1$ and all other pairs of $(-1)$-curves are disjoint. By the Hodge index theorem, $F_i = E_i + E_i'$ are all equal and it defines a fibration
$
f: X \rightarrow C;
$
(iii) $E$ is a smooth rational curve such that $E^2 = 0$ and its Galois conjugates are disjoint to each other. It defines a fibration $f:X\rightarrow C$.

Taking the $\tau$-invariant part of this equality, we obtain
\[
^\tau\overline{\mathrm{NE}}(\overline{X})^G = {}^\tau\overline{\mathrm{NE}}(\overline{X})^G_{K_X \geq0} + \sum_{L \in \mathcal F} \R_+ [L+\tau L],
\]
where $\mathcal F$ is a subset of $\mathcal E$ such that for any $L\in \mathcal F$, $[L+\tau L]$ generates an extremal ray for $^\tau\overline{\mathrm{NE}}(\overline{X})^G$. Note that the action of $\tau$ and $G$ commute because $\tau$ is defined over $k$. 

When $[L+\tau L]$ generates an extremal ray, one needs to have $(L+\tau L)^2 \leq0$. Thus we have the following possibilities: (i-a) $L$ is as in case (i), $L \cap \tau L = \emptyset$ or $\tau L = L$; (i-b) $L$ is as in case $(i)$, and we have $(L+\tau L)^2 = 0$; (ii-a) $L$ is as in case (ii) and we have $[L] = [\tau L]$; (iii-a) $L$ is as in case (iii) and we have $[L] = [\tau L]$.

In the situations of (ii-a) and (iii-a), our assertion is immediate. The case of (i-a) does not occur since we are assuming that the pair $(X, \tau)$ is minimal. In the case of (i-b), let $F = L+\tau L$. There exists the unique Galois conjugate $E'$ of $E$ such that $E \cdot \tau E' =1$, and this implies that $F$ is a nef divisor. By Riemann--Roch, we have $h^0(F) \geq 2$. Now $L$ and $\tau L$ do not move linearly, so $|F|$ is a base point free linear system with desired properties.
\end{proof}

Using the previous lemma, we can conclude the following proposition:

\begin{prop}
\label{proposition:conicfibrations}
Let $(X, \tau)$ be a minimal pair with $\mathrm{rk}\,\NS(X)^\tau >1$ and $K_X$ not nef. Then $(X, \tau)$ is one of the following:
\begin{enumerate}
\item There exists a smooth conic fibration $f: X\rightarrow C$ with a non-trivial involution $\upsilon$ of $C$ such that $f\circ \tau = \upsilon \circ f$ and $C$ is a smooth projective curve;
\item There exists a conic fibration $f: X\rightarrow C$, with singular fibers with a non-trivial involution $\upsilon$ of $C$ such that $f\circ \tau = \upsilon \circ f$ and $C$ is a smooth projective curve;
\item There exists a conic fibration $f: X \rightarrow C$ such that $f\circ \tau = f$; the involution $\tau$ induces a non-trivial involution on every smooth fiber and $C$ is a smooth projective curve.
\end{enumerate}
\end{prop}
\begin{proof}
By Lemma~\ref{lemma:bpfpencil}, we have a $\tau$-invariant base point free linear system of Iitaka dimension $1$ defined over $k$, and this defines a fibration $X \rightarrow C$ with an involution $\upsilon$ of $C$ (possibly trivial) such that $\upsilon \circ f = f \circ \tau$. The minimality implies that $f$ is a conic fibration.

Assume that $\upsilon$ is trivial. We would like to prove that $\tau$ induces a non-trivial action on every smooth fiber so that this case gives case (3). Let $D$ be the fixed locus of $\tau$. The locus $D$ contains the horizontal curve $D_h$ and the vertical curve $D_v$. Since the action of $\tau$ on the generic fiber of $f$ is non-trivial, the map $D_h \rightarrow C$ has to be a degree-two finite morphism. Let $U$ be the complement of singular fibers of $f$. Consider the following morphism:
\[
\tilde{f}: U/\langle \tau \rangle \rightarrow f(U)
\]
Since every fiber of $f : U  \to U/\langle \tau \rangle \to  f(U)$ is reduced, we conclude that $D_v \cap U$ is the empty set. 
Thus our assertionn follows
\end{proof}

When we have $\mathrm{rk} \,\NS(X)^\tau =1$, we have the following proposition instead:

\begin{prop} \label{prop:rankonenotnef}
Let $(X, \tau)$ be a minimal pair with $\mathrm{rk}\,\NS(X)^\tau =1$ and $K_X$ not nef. Then $X$ is a del Pezzo surface. Let $r$ be the positive integer such that $(r) = (K_X, \NS(X))$ where $(K_X, \NS(X))$ is the ideal of $\mathbb Z$ generated by $K_X.D$ for $D$ a divisor on $X$. 
The surface $X$ satisfies one of the following statements:
\begin{enumerate}
\item $X$ is minimal over $k$ with $\rk \Pic(X) = 1$ or
a geometrically quadric surface with $\tau$ interchanging the rulings;
\item $X$ is a non-minimal del Pezzo surface of degree $6$; more precisely $X$ is a blow up of a closed point of degree $3$ on $X'$ where $X'$ is geometrically isomorphic to $\PP^2$; we have $r = 3$;
\item $X$ is a non-minimal del Pezzo surface of degree $4$ and we have $r = 2$ or $4$;
\item $X$ is a non-minimal cubic surface and we have $r = 3$;
\item $X$ is a non-minimal del Pezzo surface of degree $2$ and $r=1$ or $2$;
\item $X$ is a non-minimal del Pezzo surface of degree $1$ and $r = 1$.
\end{enumerate}
\end{prop}
Here $\Pic(X)$ denotes classes of divisors $D\subset X$ defined over $k$.

\begin{proof}
It is easy to see that $-K_X$ is ample.
If $\mathrm{rk} \, \mathrm{Pic}(X) =1$, then $X$ is minimal over $k$. 

Suppose that $\mathrm{rk} \, \mathrm{Pic}(X) >1$. 
The action of $-\tau^*$ on $\mathrm{Pic}(X)$ is the reflection with respect to $K_X^\perp$ \cite{Wal87}*{\S 2}. 
Let $\alpha$ be the generator of the ray containing $K_X$. If $K_X$ is divisible, then $X$ is isomorphic to a quadric surface. Thus it is minimal over $k$. We may assume that $\alpha = K_X$ and $X$ is not minimal over $k$. Then we must have
\[
r \mid (K_X, K_X) \mid 2r.
\]
Indeed, there exists $L \in \mathrm{Pic}(X)$ such that $(K_X, L)=r$.
Then $(K_X, L + \tau L) = 2r$, and $L+\tau L$ is a multiple of $K_X$, proving the claim.

Assume that $(K_X, K_X) = 8$. One possibility is that $X$ contains a unique exceptional curve over the ground field. We conclude that $r = 1$. This is a contradiction. When $X$ is geometrically a quadric surface, $\tau$ must interchange two rulings.

Suppose that $(K_X, K_X) = 7$. Since $X$ is not minimal, we may have an exceptional curve defined over the ground field on $X$ or the disjoint union of two exceptional curves conjugate to each other. This means that $r = 1, 2$, which is impossible.

Assume that $(K_X, K_X) = 6$. Since $X$ is not minimal, the possibilities of $r$ are $1, 2$, or $3$. Hence we conclude that $r =3$, and there exists a birational morphism $f$ to $X'$ such that $X'$ is isomorphic to $\PP^2$ over the algebraic closure. The morphism $f$ is a blow up of a closed point of degree $3$ on $X'$. 

If $(K_X, K_X) = 5$, then the possibilities of $r$ are $1, 2, 3$, or $4$. This is impossible.

Other cases are similar.
\end{proof}

\subsection*{The classification over the reals}
Over $\R$, we state more precise statements:

\begin{prop}
\label{proposition:conicfibrations over real}
Let $(X, \tau)$ be a minimal pair over $\R$ with $\mathrm{rk} \, \NS(X)^\tau >1$ and $K_X$ not nef. Then $(X, \tau)$ is one of the following:
\begin{enumerate}
\item There exists a smooth conic fibration $f: X\rightarrow C$ to a smooth projective curve $C$ with a non-trivial involution $\upsilon$ of $C$ such that $f\circ \tau = \upsilon \circ f$;
\item There exists a conic fibration $f: X\rightarrow C$ to a smooth projective curve $C$, with singular fibers, with a non-trivial involution $\upsilon$ of $C$ such that $f\circ \tau = \upsilon \circ f$. Let $F_0$ be a singular fiber for $f$. Then we have either
\begin{enumerate}
\item the fiber $F_0$ has the form of $E + E'$ where $E$ is a $(-1)$-curve and $E'$ is its Galois conjugate with $E\cdot E' =1$, or ;
\item the fiber $F_0$ has the form of $E + \tau E'$ where $E$ be an $(-1)$-curve in $F_0$, then its Galois conjugate $E'$ is disjoint with $E$ and we have $E\cdot \tau E' = 1$;
\end{enumerate}
\item There exists a conic fibration $f: X \rightarrow C$ to a smooth projective curve $C$ such that $f\circ \tau = f$; the involution $\tau$ induces a non-trivial involution on smooth fibers. Let $F_0$ be a singular fiber. Then we have
\begin{enumerate}
\item the fiber $F_0$ has the form of $E + E'$ where $E$ is a $(-1)$-curve and $E'$ is its Galois conjugate with $E\cdot E' =1$, or ;
\item the fiber $F_0$ has the form of $E + \tau E$ where $E$ be an $(-1)$-curve in $F_0$, and we have $E\cdot \tau E = 1$;
\end{enumerate}
\end{enumerate}
\end{prop}

\begin{proof}
This follows from Proposition~\ref{proposition:conicfibrations} and its proof. There is one thing we need to show: In case (2), any singular fiber $F_0$ is not of the form of $E+\tau E$ where $E$ is a $(-1)$-curve with $E \cdot \tau E = 1$. 

Suppose that it takes this form. Then the involution $\tau$ fixes the intersection point $\{p\} = E\cap \tau E$. Then one can find analytic coordinates $\{x, y\}$ such that $p$ is given by $x=y=0$ and the action $\tau$ at $p$ is given by $(x, y) \mapsto (y, x)$. Thus we can find an analytic hypersurface at $p$ which $\tau$ acts trivially. In other words, there exists a principal ideal $(f) \subset \widehat{\mathcal O}_p$ in the completion of the local ring at $p$ such that $(f)$ is preserved by $\tau^*$ and the action on $\widehat{\mathcal O}_p/(f)$ is trivial. We claim that $(f) \cap \mathcal O_p$ is non-trivial. If not, then we have the injection map $\mathcal O_p \hookrightarrow \widehat{\mathcal O}_p/(f)$, so we conclude that the action on $X$ itself is trivial and it contradicts with the fact that $\tau$ is non-trivial. This discussion shows that there exists a horizontal fixed curve passing through $p$ where the action of $\tau$ is trivial. It follows that $\upsilon$ must be trivial. This is a contradiction with our assumption.
\end{proof}

Here we show that case (1) of Proposition~\ref{proposition:conicfibrations over real} occurs.

\begin{example}
Let $X = C_1 \times C_2$ where $C_i$ is a smooth conic.
We consider an involution on $C_2$ which induces an involution on $X$ by acting on $C_1$ factor trivially.
Then the projection $X \to C_2$ is an example for the case (1).
\end{example}

Here we show that the case (2) of Proposition~\ref{proposition:conicfibrations over real} occurs.

\begin{example}
\label{example:(2-a)(2-b)}
Let $f$ and $g$ be homogeneous polynomials of degree $2$ in three variables $x_0, x_1, x_2$ over $\R$ such that $f$ and $g$ define two smooth plane conics meeting transversely. We consider the following surface:
\[
X \subset \PP^2 \times \PP^1 : s^2f(x_0, x_1, x_2) + t^2g(x_0, x_1, x_2) = 0
\]
This defines a smooth surface, and the projection $\pi$ to $\PP^1$ gives us a conic bundle structure. We consider the following involution on $\PP^1$
\[
\tau : (s: t) \mapsto (s: -t).
\]
This defines an involution $\tau$ on $X$ by mapping as
\[
    ((x_0:x_1:x_2), (s:t))\mapsto ((x_0 : x_1 : x_2), (s:-t)).
    \]

Suppose that $f$ is given by $x_0^2+x_1^2+x_2^2$. Without loss of generality, we may assume that $g$ is given by $-a_0x_0^2 - a_1x_1^2 -a_2x_2^2$ where $a_0, a_1, a_2$ are mutually distinct. Suppose that $a_0 > a_1> 0 > a_2$. Then a pair $(X, \tau)$ is not relatively minimal over $\R$. 

Indeed, $\pi$ has six geometric singular fibers at $(s:t) = (\pm \sqrt{a_0} : 1), (\pm \sqrt{a_1} : 1), (\pm \sqrt{a_2}:1)$, and the fiber at $(\sqrt{a_1}:1)$ consists of two components and we denote an exceptional curve on this fiber by $E_1$. Then we can contract $E_1 + \tau E_1$ which is Galois invariant and we obtain a pair $(X', \tau')$ with a conic bundle structure $\pi'$. We claim that this pair $(X', \tau')$ is minimal over $\R$. 

First note that $\mathrm{rk} \, \Pic (\overline{X}') = 6$ and its generators over $\mathbb Q$ are
\[
E_0, \tau E_0, E_2, \tau E_2, F, K_{X'}
\]
where $E_0$ is an exceptional curve defined over $\mathbb C$ in the fiber at $(\sqrt{a_0} : 1)$, $E_2$ is an exceptional curve defined over $\mathbb C$ in the fiber at $(\sqrt{a_2} : 1)$, and $F$ is a general fiber for $f$. We denote the Galois conjugates of $E_0$ and $E_2$ by $E_0'$ and $E_2'$ respectively. Then we have $E_0+E_0' \sim F$ and $E_2+\tau E_2' \sim F$. Thus we have $\rk \, \Pic(X')^{\tau'} = 2$ and it is generated by $F$ and $K_{X'}$ over $\mathbb Q$. The involution $\tau'$ acts on this space trivially. 

To see the minimality of $(X', \tau')$, we need to understand extremal rays of $\overline{\mathrm{NE}}(X')^G$. One ray is generated by $F$. Another ray is generated by a sum of disjoint two $(-2)$-curves. Indeed, the anticanonical divisor $-K_X$ is a big and nef divisor contracting four disjoint $(-2)$-curves which are sections for $\pi$. It follows that $-K_{X'}$ is also a big and nef divisor contracting two disjoint $(-2)$-curves. 
This implies that there is no equivariant birational contraction contracting the disjoint union of $(-1)$-curves. Indeed, if such a contraction exists, then the sum of $(-1)$-curves has to span a ray of $\overline{\mathrm{NE}}(X')^G$.
Thus our assertion follows.

Now note that the fibers at $(\pm \sqrt{a_0} : 1)$ are type (2-a) and the fibers at $(\pm \sqrt{a_2} : 1)$ are type (2-b). So this is an example for both (2-a) and (2-b).
\end{example}

\begin{example}
    Let $X$ be the surface considered in Example~\ref{example:(2-a)(2-b)}.
    Let $\zeta$ be the involution on $X$ defined by mapping as
    \[
    ((x_0:x_1:x_2), (s:t))\mapsto ((x_0 : x_1 : -x_2), (s:t)).
    \]
    Assume that we have $a_0 > a_1 > a_2 > 0$. Then we have $\rk \Pic(X)^{\zeta} = 2$. As before, one can show that $(X, \zeta)$ is minimal. Fibers at $(\pm\sqrt{a_0} : 1)$ and $(\pm\sqrt{a_2}:1)$ are witnesses to (3-a) and fibers at $(\pm\sqrt{a_1}:1)$ are witnesses to (3-b). 
\end{example}

\begin{prop}
\label{proposition:delpezzo}
Retain the notation of Proposition~\ref{prop:rankonenotnef}. Then $X$ is a del Pezzo surface classified as follows:
\begin{enumerate}
\item $X$ is minimal over $\R$ with $\rk \Pic(X) = 1$ or
a geometrically quadric surface with $\tau$ interchanging the rulings;
\item $X$ is a non-minimal del Pezzo surface of degree $4$ and we have $r = 2$;
\item $X$ is a non-minimal del Pezzo surface of degree $2$ and $r=1$ or $2$;
\item $X$ is a non-minimal del Pezzo surface of degree $1$ and $r = 1$.
\end{enumerate}
\end{prop}
\begin{proof}
A proof is the same as before.
Note that over $\R$, there are only closed points of degree $1$ and $2$. Thus the case $r=K_X^2=4$ -- a non-split quadric surface blown up along a degree-four closed point -- does not arise over $\R$.
\end{proof}

\begin{example}
    Let us exhibit examples for each case in Proposition~\ref{proposition:delpezzo}. An involution on $\PP^2$ or the involution of a degree-two morphism 
    $X \rightarrow \PP^2$ branched
    along a conic are witnesses of (1)
    For a non-minimal del Pezzo surface of degree $2$, one can consider the Geiser involution. For a non-minimal del Pezzo surface of degree $1$ one can consider the Bertini involution.

    Finally for a non-minimal del Pezzo surface of degree $4$, let $Q_{3, 1}$ be a quadric surface of Picard rank $1$. We pick a smooth member $C \in |-K_{Q_{3,1}}|$ and let $\pi : X \to Q_{3,1}$ be the double cover of $Q_{3,1}$ ramified along $C$. Let $\tau$ be the involution associated to this covering. Then we have $\rk \Pic (X)^\tau = 1$. Let $\ell_1 \subset Q_{3,1}$ be a line which is tangent to $C$. (This cannot be defined over $\R$ as the Picard rank of $Q_{3,1}$ is $1$.) Let $\ell_2$ be the Galois conjugate of $\ell_1$. The pullback of $\ell_1$ is denoted by $E_1 + \tau E_1$ and the pullback of $\ell_2$ is $E_2 + \tau E_2$ where $E_2$ is a conjugate of $E_1$. If necessary, after replacing $X \to Q_{3,1}$ by its twisted form, we may assume that $E_1 \cap E_2 = \emptyset$.
    Indeed, if $E_1\cap E_2\neq \emptyset$, then after twisting the conjugate of $E_1$ becomes $\tau E_2$ which satisfies $E_1\cap \tau E_2 = \emptyset$. Thus one can contract $E_1 + E_2$  and it shows that $X$ is non-minimal and $r = 2$.

    Note that assuming $X(\R) \neq \emptyset$, $X$ is rational as any del Pezzo surface of degree $\geq 5$ with a rational point is rational over the ground field.
    Indeed, since $X$ is not minimal, $X$ is birational to a del Pezzo surface of degree $\geq 5$ with a real point. Thus $X$ is rational.
    In particular this involution has been studied in \cite{CMYZ24}*{Section 9}.
\end{example}

\section{Degree-Two Unirational Parametrizations}
\label{sec:Kollar}

We assume that our ground field is $\R$.
Let $S$ be a minimal conic bundle with $2m \, (m\geq 2)$ singular fibers defined over $\R$.
Then it is well-known that $S$ is non-rational over $\R$. (See, e.g., \cite{Kollar}.)

Here we recall the following theorem of Koll\'ar and Szab\'o:

\begin{thm}[{\cite{KS03}}]
\label{thm:bisection}
Let $C$ be a smooth projective geometrically integral curve of genus $g$ defined over $\R$. Let $\pi : W \to C$ be a dominant morphism from a smooth projective variety $W$ defined over $\R$ such that the generic fiber $W_\eta$ is geometrically rationally connected. 
Suppose that there is a real closed point $w \in W(\R)$ such that the fiber $W_{\pi(w)}$ containing $w$ is smooth.
Then there exists a bisection of $\pi$ of geometric genus $2g$ defined over $\R$.
\end{thm}

\begin{proof}
Since $W_\eta$ is geometrically rationally connected, it follows from \cite{GHS03} and \cite{HT06} that there is a section $\Sigma \subset W$ defined over $\mathbb C$ such that $\Sigma$ contains $w$ and its normal bundle is globally generated and its $H^1$ vanishes. Let $\overline{\Sigma}$ be the complex conjugate of $\Sigma$ and we consider the following prestable curve of genus $2g$
\[
C' = \Sigma \cup \PP^1 \cup \overline{\Sigma}.
\]
Here $\Sigma$ is glued to $\PP^1$ at $w \in \Sigma$ and $p \in \PP^1(\mathbb C) \setminus \PP^1(\R)$, and $\overline{\Sigma}$ is glued to $\PP^1$ at $w \in \overline{\Sigma}$ and the complex conjugate $\overline{p} \in \PP^1$.
Then we consider the following stable map of genus $2g$
\[
f : C' \to W \times \PP^1,
\]
mapping $\Sigma$ to $\Sigma \times \{p\}$, $\PP^1$ to $\{w\} \times \PP^1$, and $\overline{\Sigma}$ to $\overline{\Sigma}\times \{\overline{p}\}$. This corresponds to a smooth real point of the moduli space of stable maps.

Then it follows from \cite{GHS03}*{Lemma 2.6} that one can smooth this stable map to a stable map of genus $2g$ $f: C' \to W \times \PP^1$ defined over $\R$. The resulting stable map $f' : C' \to W \times \PP^1 \to W$ is a bisection of $\pi$. Thus our assertion follows.
\end{proof}

The following corollary was first obtained in \cite{Isk67}*{Corollary 4.4}, and the above theorem also proves it:

\begin{cor}
\label{cor:conicbundleunirational}
Let $\pi : S \to \PP^1$ be a conic bundle defined over $\R$ such that $S$ is smooth and $S(\R)$ is non-empty. Then $S$ admits a degree $2$ unirational parametrization.
\end{cor}
\begin{proof}
It follows from Theorem~\ref{thm:bisection} that $\pi$ admits a bisection of genus $0$. The base change of $\pi$ by this bisection will give us a degree $2$ unirational parametrization of $S$.
\end{proof}

\section{The proof of the main theorem}
\label{sec:main}

The goal of this section is to prove Theorem~\ref{thm:intromain}.
To this end, we will show that if $(X, \tau)$ is a pair of a smooth real rational projective surface with an involution $\tau$ such that $(X, \tau)$ is minimal, then $X/\tau$ is either rational or birational to a minimal conic bundle using the classification of such pairs in \cite{CMYZ24}.

Let us recall the classification of birational involutions of the real projective plane provided in \cite{CMYZ24}.
For a birational involution $\iota \in \Bir(\PP^2)=\Bir(\PP^2_\R)$ and its regularization $\tau \in \Aut(X)$ on some smooth $\R$-rational surface $X$, denote by $F(\tau)$ the union of all geometrically irrational real curves in the surface $X$ that are pointwise fixed by the involution $\tau$.
The classification is as follows.

\begin{thm}[{\cite{CMYZ24}*{Main Theorem}}]
\label{thm:classification of involutions over real}
Let $\iota$ be an involution in $\Bir(\PP^2)$.
Then $\iota$ admits a regularisation $\tau \in \Aut(X)$ on a smooth real rational projective surface $X$ such that a pair $(X, \tau)$ with $G = \langle \tau \rangle$ belongs to one of the following classes:

\begin{itemize}

\item[$(\frakL)$] The surface $X$ is $\PP^2$ and $\tau$ is the linear involution $(x : y : z) \to (x : y : -z)$ with $F(\tau)=\emptyset$.

\item[$(\frakQ)$] The surface $X$ is a quadric surface in $\PP^3$, and $\tau$ is an involution with $F(\tau)=\emptyset$ and $X(\R)^\tau = \emptyset$.

\item[$(\frakT_{4n})$] The surface $X$ admits a $G$-equivariant morphism $X \to \PP^1$ that is a conic bundle with $4n \geq 4$ singular fibers and $\Pic(X)^G \simeq \Z^2$, and the involution $\tau$ is a $0$-twisted Trepalin involution with $F(\tau)=\emptyset$.

\item[$(\frakT'_{4n+2})$] The surface $X$ admits a $G$-equivariant morphism $X \to \PP^1$ that is a conic bundle with $4n+2 \geq 6$ singular fibers and $\Pic(X)^G \simeq \Z^2$, and the involution $\tau$ is a $1$-twisted Trepalin involution with $F(\tau)=\emptyset$.

\item[$(\frakT''_{4n})$] The surface $X$ admits a $G$-equivariant morphism $X \to \PP^1$ that is a conic bundle with $4n \geq 4$ singular fibers and $\Pic(X)^G \simeq \Z^2$, and the involution $\tau$ is a $2$-twisted Trepalin involution with $F(\tau)=\emptyset$. 

\item[$(\frakB_4)$] The surface $X$ is a del Pezzo surface of degree $1$ with $\Pic(X)^G \simeq \Z$, and the involution $\tau$ is the Bertini involution of the surface $X$ such that $F(\tau)$ is a non-hyperelliptic curve of genus $4$.

\item[$(\frakG_3)$] The surface $X$ is a del Pezzo surface of degree $2$ with $\Pic(X)^G \simeq \Z$, and the involution $\tau$ is the Geiser involution of the surface $X$ such that $F(\tau)$ is a non-hyperelliptic curve of genus $3$.

\item[$(\frakK_1)$] The surface $X$ is a del Pezzo surface of degree $2$ with $\Pic(X)^G \simeq \Z$, and the involution $\tau$ is the Kowalevskaya involution of the surface $X$ such that $F(\tau)$ is a genus $1$ curve.
\item[$(\frakd\frakJ_g)$] The surface $X$ admits a $G$-equivariant morphism $X \to \PP^1$ that is a $G$-exceptional conic bundle with $2g+2$ singular fibers and $\Pic(X)^G \simeq \Z^2$, and the involution $\tau$ is a de Jonqui\`eres involution such that $F(\tau)$ is a hyperelliptic curve $C$ of genus $g \geq 1$.

\item[($\frakI_g$)] The surface $X$ admits a $G$-equivariant morphism $X \to \PP^1$ that is a non-$G$-exceptional conic bundle with $2g+2$ singular fibers and $\Pic(X)^G \simeq \Z^2$, and the involution $\tau$ is a $0$-twisted Iskovskikh involution such that $F(\tau)$ is a hyperelliptic curve $C$ of genus $g \geq 1$.

\item[$(\frakI'_g)$] The surface $X$ admits a $G$-equivariant morphism $X \to \PP^1$ that is a non-$G$-exceptional conic bundle with $2g+3$ singular fibers and $\Pic(X)^G \simeq \Z^2$, and the involution $\tau$ is a $1$-twisted Iskovskikh involution such that $F(\tau)$ is a hyperelliptic curve $C$ of genus $g \geq 1$.

\item[$(\frakI''_g)$] The surface $X$ admits a $G$-equivariant morphism $X \to \PP^1$ that is a non-$G$-exceptional conic bundle with $2g+4$ singular fibers and $\Pic(X)^G \simeq \Z^2$, and the involution $\tau$ is a $2$-twisted Iskovskikh involution such that $F(\tau)$ is a hyperelliptic curve $C$ of genus $g \geq 1$. 

\end{itemize}

\end{thm}

In the following, we consider the quotient $X / \tau$ for a pair of a smooth real projective rational surface $X$ and an involution $\tau$ on $X$ such that $(X, \tau)$ is one of pairs listed in Theorem~\ref{thm:classification of involutions over real}.

To begin with, the quotient $\PP^2/\tau$ of $\PP^2$ by the linear involution $\tau$ which acts by $(x : y : z) \to (x : y : -z)$ is rational because we can take $y/x, (z/x)^2$ as a transcendental basis of the function field of $\PP^2/\tau$ over $\R$.
This is the case of the class $\frakL$ of Theorem~\ref{thm:classification of involutions over real}.

For the class $\frakQ$, we have the following proposition:

\begin{lemma}
Let $X$ be a rational quadric surface in $\PP^3$ and let $\tau$ be an involution on $X$.
Then the quotient $X/\tau$ is rational.
\end{lemma}
\begin{proof}
According to \cite{CMYZ24}*{Proposition 3.1}, there are only three possibilities:
\begin{enumerate}
\item $\tau$ is conjugate to a linear involution of $\PP^2$ in $\Bir(\PP^2)$.
\item $X \simeq \{x^2+y^2+z^2=w^2\}$ and $\tau$ is the antipodal involution $(x:y:z:w) \mapsto (x:y:z:-w)$.
\item $X \simeq \PP^1 \times \PP^1$ and $\tau$ acts by $((x:y),(s:t)) \mapsto ((x:y),(t:-s))$.
\end{enumerate}
In case (1), the quotient $X/\tau$ is rational because the quotient of $\PP^2$ by a finite group of automorphisms is rational (see, e.g.\cite{Tre14}).
In case (2), the function field of $X/\tau$ is $\R(Y,Z,W)^{(Y,Z,W)\mapsto (Y,Z,-W)}= \R(Y,Z,W^2)$ with $1+Y^2+Z^2=W^2$, where $Y=y/x, Z=z/x, W=w/x$.
Thus, $X/\tau$ is rational.
In the case (3), $\tau$ acts trivially on the first factor of $X\simeq \PP^1\times \PP^1$.
The quotient of the second factor has a real point so it is rational over $\mathbb R$.
So, $X/\tau$ is rational.
\end{proof}

For the Bertini involution, we have the following lemma:

\begin{lemma}
Let $X$ be a rational del Pezzo surface of degree $1$ and let $\tau$ be the Bertini involution on $X$.
Then the quotient $X / \tau$ is rational.
\end{lemma}
\begin{proof}
The Bertini involution $\tau$ on $X$ is the Galois involution of the double covering $X \to Q \subset \PP^3$ defined by the divisor $-2K_X$, where $Q$ is a geometrically irreducible quadric cone in $\PP^3$.
Thus we have $X / \tau \simeq Q$.
Since $X$ is rational, $Q$ has plenty of smooth $\R$-points.
Hence, our assertion follows.
\end{proof}

Next we handle the Geiser involution:

\begin{lemma}
Let $X$ be a del Pezzo surface of degree $2$ and let $\tau$ be the Geiser involution on $X$.
Then the quotient $X/\tau$ is rational.
\end{lemma}
\begin{proof}
The Geiser involution $\tau$ on $X$ is the Galois involution of the double covering $X \to \PP^2$ defined by the anticanonical divisor $-K_X$.
Thus we have $X / \tau \simeq \PP^2$.
\end{proof}

Finally we discuss the Kowalevskaya involution:

\begin{lemma}[{\cite{CMYZ24}*{Remark 6.6}}]
Let $X$ be a rational del Pezzo surface of degree $2$ and let $\tau$ be the Kowalevskaya involution on $X$. 
Then the quotient $X/\tau$ is birationally equivalent to a minimal Iskovskikh surface (i.e. a minimal exceptional conic bundle with $4$ singular fibers).
\end{lemma}

\begin{proof}
Let us recall that the Kowalevskaya involution $\tau$ acts on a del Pezzo surface $X$ in $\PP(1,1,1,2)$ given by
\[
w^2=-(y^2+ax^2+bxz+cz^2)^2\pm xz(x-z)(x-sz)
\]
for some real numbers $a, b, c, s$ with $a>0, s>1$, and the involution $\tau$ is given by
\[
(x:y:z:w) \mapsto (x:-y:z:w).
\]
Thus the quotient $S=X/\tau$ is a hypersurface in $\PP(1,1,2,2)$ that is given by
\[
w^2=-(u+ax^2+bxz+cz^2)^2\pm xz(x-z)(x-sz).
\]
The surface $S$ is a del Pezzo surface of degree 4 that is singular at the points $(0:0:1:i)$ and $(0:0:1:-i)$ and its singularities at these points are ordinary double points.
Blowing up the conjugated points $(0:0:1:i)$ and $(0:0:1:-i)$, we get a conic bundle $\eta \colon \widetilde{S} \to \PP^1$, which is defined over $\R$.
A conic bundle $\eta$ has exactly four geometrically singular fibers.
Namely, they are the preimages of the curves in $S$ that are cut out on $Y$ by the equations $x=0, z=0, x=z, x=sz$.
These fibers are conics in $\PP^2$ that are irreducible over $\R$.
This implies that $\rk \Pic(\widetilde{S})=2$.
\end{proof}

Finally we record the following theorem in the case of equivariant conic bundles:

\begin{thm}[{\cite{Tre16}*{Theorem 4.1}}]
Let $G$ be a finite group.
Let $X$ be a $G$-smooth projective surface that admits a $G$-equivariant conic bundle structure.
Then $X/G$ is birational to a relatively minimal conic bundle.
\end{thm}

\begin{proof}[Proof of Theorem~\ref{thm:intromain}]
    Results in this section show that if $(X, \tau)$ is a pair of a smooth real rational projective surface with an involution $\tau$, then $X/\tau$ is either rational or birational to a minimal conic bundle. This shows that when a real geometrically rational surface $S$ admits a degree $2$ unirational parametrization, $S$ is either rational or birational to a minimal conic bundle. Conversely when $S$ is rational, it is clear that it admits a degree $2$ unirational parametrization. For a conic bundle $S$, it follows from Corollary~\ref{cor:conicbundleunirational} that $S$ admits a degree $2$ unirational parametrization. Thus our assertion follows.
\end{proof}

\begin{proof}[Proof of Corollary~\ref{cor:degree12}]
    This follows from Theorems~\ref{thm:intromain} and \ref{thm:kollar}.
\end{proof}

\begin{example}
    Let us explain how to construct minimal del Pezzo surfaces of degree $2$. Let $S'$ be a blow up of $7$ general real points on $\PP^2$. This is a split del Pezzo surface of degree $2$. The anticanonical linear system defines a double cover $S' \to \PP^2$ ramified along a smooth quartic curve $C$. Since $S'$ is split, all $(-1)$-curves are defined over $\R$. This means  that all 28 bitangent lines to $C$ are defined over $\R$, and each pullback consists of two $(-1)$-curves defined over $\R$. Let $S \to  \PP^2$ be the twisted form of $S' \to \PP^2$ by the Geiser involution. Then the pullback of each bitangent line consists of two $(-1)$-curves conjugated to each other, and we conclude that $\rk \Pic(S) = 1$. Thus this is a minimal del Pezzo surface.
    
    Conversely if we have a minimal del Pezzo surface $S \to \PP^2$ of degree $2$ and Picard rank $1$, then for any $(-1)$-curve $E$ and its complex conjugate $E'$, we have $E + E' \sim -K_S$. In particular all bitangent lines are defined over $\R$, each pullback consists of two $(-1)$-curves conjugated to each other.
    Such a del Pezzo surface admits a unirational parametrization of degree $24$ by \cite{Manin}*{Theorem 29.9}. However, our theorem shows that it does not admit a degree-two real unirational parametrization.

\end{example}
\begin{example}
    Next let us explain how to construct minimal del Pezzo surfaces of degree $1$. Let $S'$ be a blow up of $8$ general real points on $\PP^2$. This is a split del Pezzo surface of degree $1$. The linear system $|-2K_S|$ defines a double cover $S' \to Q$ to a quadric cone $Q$ ramified along a smooth sextic curve $C$ and the vertex of $Q$. Since $S'$ is split, all $(-1)$-curves are defined over $\R$. This means  that all tritangent conics to $C$ are defined over $\R$, and each pullback consists of two $(-1)$-curves defined over $\R$. Let $S \to Q$ be the twisted form of $S' \to Q$ by the Bertini involution. Then the pullback of each tritangent conic consists of two $(-1)$-curves conjugated to each other, and we conclude that $\rk \Pic(S) = 1$. Thus this is minimal. 
    
    Conversely if we have a minimal del Pezzo surface $S \to Q$ of degree $1$ and Picard rank $1$, then for any $(-1)$-curve $E$ and its complex conjugate $E'$, we have $E + E' \sim -2K_S$. In particular all tritangent conics are defined over $\R$, and each pullback consists of two $(-1)$-curves conjugated to each other. It is a major open problem whether $S$ admits a unirational parametrization.
\end{example}

Finally to prove Corollary~\ref{cor:degreegeq3}, we prepare the following lemma:
\begin{lemma}
\label{lemma:degree3}
    Let $S$ be a real del Pezzo surface of degree $3 \leq d \leq 7$. Then either $S$ is not minimal or $S$ admits a conic bundle structure over $\PP^1$.
\end{lemma}

\begin{proof}
    Assume that $S$ is minimal. Let $E_1$ be a line on $S$ and $E_2$ be its conjugate. Then since the degree of $S$ is greater than or equal to $3$, we must have $E_1.E_2 = 0$ or $1$. However, the minimality implies we must have $E_1.E_2 = 1$. This means that $E_1 + E_2$ defines a conic bundle structure over $\PP^1$. Thus our assertion follows.
\end{proof}

\begin{proof}[Proof of Corollary~\ref{cor:degreegeq3}]
    It follows from Lemma~\ref{lemma:degree3} that $S$ is birational to a del Pezzo surface $S'$ of degree $\geq 3$ which is minimal. If the degree of $S'$ is less than or equal to $7$, then this is a minimal conic bundle. Thus Theorem~\ref{thm:intromain} implies that $S'$ admits a degree-two real unirational parametrization. When the degree $S' \geq 8$, the existence of a real point on $S$ implies that $S'$ is rational. Thus our assertion follows. 
\end{proof}


\end{document}